\documentclass{article}
\usepackage{amsmath,amsthm,amsfonts}
\usepackage{amssymb}
\usepackage{graphicx}

\usepackage{color}

\newtheorem{subthm}[subsubsection]{Theorem}

\newtheorem{subprop}[subsubsection]{Proposition}
\newtheorem{subdef}[subsection]{Definition}
\newcommand{\bracket}[1]{ \left\{\begin{array}{l} #1 \end{array} \right.}

\newcommand{\ind}[1]{1_{#1}}
\newcommand{\Nats}{\mathbb{N}}
\newcommand{\Reals}{\mathbb{R}}
\newcommand{\espc}{\mathbb{E}}
\newcommand{\eps}{\epsilon}
\newcommand{\esp}{\mathbb{E}}
\newcommand{\calF}{\mathcal F}

\newcommand{\E}{{\mathbb{E}}}

\newtheorem{lem}[subsubsection]{Lemma}

\newcommand{\calU}{\mathcal U}
\newcommand{\eref}[1]{Equation~(\ref{#1})}
\newcommand{\sref}[1]{Section~\ref{#1}}
\def\P{\mathbb{P}}

\newcommand{\calA}{\mathcal A}

\newcommand{\calKK}{\mathcal B}
\newcommand{\calS}{\mathcal S}
\newcommand{\calX}{\mathcal X}
\def\be{\begin{equation}}
\def\ee{\end{equation}}
\def\ben{\[}
\def\een{\]}
\def\bearn{\begin{eqnarray*}}
\def\eearn{\end{eqnarray*}}
\def\bear{\begin{eqnarray}}
\def\eear{\end{eqnarray}}
\def\barr{\begin{array}}
\def\earr{\end{array}}
\begin{document}
%%%\title{Mean Field Markov Games}
\title{Mean Field Asymptotics of Markov Decision Evolutionary Games and Teams
\thanks{
This work was partially supported by the INRIA ARC Program:
Populations, Game Theory, and Evolution (POPEYE) and by an EPFL PhD
internship grant. } \thanks{This paper has been presented at the first international conference on Game Theory for Networks, Gamenets 2009, Istanbul, Turkey, \cite{gamenets}. }}

\author{H. Tembine, J.-Y. Le Boudec, R. El-Azouzi, E.
Altman
}
\date{First version July 2008. This version March 2010.}
\maketitle \thispagestyle{empty} \pagestyle{empty}

\begin{abstract}
We introduce Mean Field Markov games with $N$ players, in
which each individual in a large population interacts with other
randomly selected players. The states and actions of each player in
an interaction together determine the instantaneous payoff for all
involved players. They also determine the transition probabilities
to move to the next state. Each individual wishes to maximize the
total expected discounted payoff over an infinite horizon. We
provide a rigorous derivation of the asymptotic behavior of this
system as the size of the population grows to infinity. Under indistinguishability  per type assumption, we show that
under any Markov strategy, the random process consisting of one
specific player and the remaining population converges weakly to a
jump process driven by the solution of a system of differential
equations. We characterize the solutions to the team and to the game
problems at the limit of infinite population and use these to
construct near optimal strategies for the case of a finite, but
large, number of players. We show that the large population
asymptotic of the microscopic model is equivalent to a (macroscopic)
mean field stochastic game in which a local interaction is
described by a single player against a population profile (the mean field limit). We
illustrate our model to derive the equations for a dynamic
evolutionary Hawk and Dove game with energy level.
\end{abstract}

%\noindent{\bf Key words:} Mean Field Interaction Model, Evolutionary
%game,, Markov Decision Precess, Chain, Game Theory, Howk-Dove game.
%%\end{IEEEkeywords}
\section{Introduction}

We consider a large population of players in which frequent interactions occur between small numbers of chosen individuals.
Each interaction in which a player is involved
can be described as one stage of a dynamic game.
The state and actions of the players at each stage
determine an immediate payoff (also called {\it fitness} in behavioral ecology) for each player
as well as the transition probabilities of a controlled
Markov chain associated with each player. Each player
wishes to maximize its expected fitness averaged over time.

This model
extends the basic evolutionary games by introducing a controlled state
that characterizes each player. The stochastic dynamic games at each
interaction replace the matrix games, and the objective of
maximizing the expected long-term payoff over an infinite time horizon replaces
the objective of maximizing the outcome of a matrix game.
Instead of a choice of a (possibly mixed) action, a player is now
faced with the choice of
decision rules (called strategies) that determine
what actions should be chosen at a given interaction for given present
and past observations.

This model with a finite number of players, called a mean field
interaction model, is in general difficult to analyze because of the
huge state space required to describe  the sate of all players. Then,
taking the asymptotics as the number of players grows to infinity,
the whole behavior of the population is replaced by  a deterministic
 limit that represents the system's state,
which is fraction of the population
at each individual state that use a given action.

In this paper we study the asymptotic {\it dynamic} behavior
of the system in which the population profile evolves in time.
For large $N,$
under mild assumptions (see Section~\ref{secmain}), the mean field converges to a deterministic measure that satisfies a non-linear ordinary differential equation for under any stationary strategy.
We show that the mean field interaction is  asymptotically equivalent to a Markov decision evolutionary game.
When the rest of the population uses a fixed strategy $u,$
any given player sees an equivalent game against
a collective of players whose state evolves according to
the ordinary differential equation (ODE)
which we explicitly compute. In addition to providing
the exact limiting asymptotic, the ODE approach provides
tight approximations for fixed large N.
The mean field asymptotic calculations for large $N$ for given choices
of strategies allows us to compute the equilibrium of the game
in the asymptotic regime.

\subsection{ Related Work}
 Mean field interaction models have already been used in
standard evolutionary games in a completely different context:
that of evolutionary game dynamics (such as replicator dynamics)
see e.g. \cite{tanabep} and references therein.
The paradigm there has been to associate relative growth rate to
actions according to the fitness they achieved,
then study the asymptotic trajectories of the state of the system,
i.e. the fraction of users that adopt the different actions.
{\it Non-atomic} Markov Decision Games have been studied in \cite{jovanovic} and
 applied  in \cite{oblivious} to firm idiosyncratic random shocks using decentralized strategies. They  proposed
the notion of {\it oblivious equilibria} via a mean field approximation. Extension to unbounded cost function can be found in \cite{ramesh}. Applications to cellular communications can found in \cite{valuetools}.
%Transition from discrete to continuous time analysis of Markov Decision Processes can be found in \cite{gast}.
% Closely related analysis of  mean field games in continuous time can be found in \cite{lions,huang,gomes}.

 Most of these approaches considered the case where the payoff of a player depends on the states of the other players but not explicitly  on the  actions of the others. In this paper, the payoff depends explicitly on both states and actions of the other players.

 %Also, it is implicitly assumed that the othe

 %In \cite{gamenets}, we introduced {\it Mean field asymptotics of Markov Decision

\subsection{ Structure}
The remainder of this paper is organized as follows. In next section we present the model assumptions and notations.  In Section~\ref{secmain}
  we present some convergence  results of the ODE in the random number of interacting players.   In Section~\ref{hd} a resource competition between animals with two types of behaviors and several states is presented.  All the sketch of proofs are given in Appendix. Section~\ref{concl}
concludes the paper.

\section{Model description}
\subsection{Mean Field Markov Process With $N$ Players}
\label{sec-model}We consider the following model, which we call
{\it Mean Field Markov  Game} with $N$ players.
%\begin{itemize}

$\bullet$ %\item
There are $N\in \mathbb{N}$ players.

$\bullet$ %\item
Each player has its own  state. A state has two components: the
\emph{type} of the player and the \emph{internal state}. The
type is a constant during the game.  The state of player $j$ at
time $t$ is denoted by $X^N_j(t)=(\theta_j,S^N_j(t))$ where $\theta_j$ is the type. The set
of possible states
 $\calX=\{1,\ldots,\Theta\}\times \calS$ is finite.

$\bullet$ %\item
 Time is discrete, taking values in
$\frac{\mathbb{N}}{N}:=\{0,\frac{1}{N},\frac{2}{N},\ldots\}.$

$\bullet$ %\item
The {\it global detailed description} of the system at time $t$
is $X^N(t)=(X^N_1(t),\ldots,X^N_N(t))$.

Define $ M^{N}(t)$ to be the current
 population profile
  i.e
$
M^{N}_{x}(t)=\frac{1}{N}\sum_{j=1}^N
\ind{\{X^N_{j}(t)=x\}}.
$
At each time $t,$ $M^N(t)$ is in the finite set
$\{0,\frac{1}{N},\frac{2}{N},\ldots,1\}^{\sharp \calX},$ and
$M^N_{\theta,s}(t)$ is the fraction of players who belong to
population of type $\theta$ (also called subpopulation
$\theta$) and have internal state $s.$ Also let
$
\bar{M}_{\theta}^N=N \sum_{s \in \calS}^N M^N_{\theta,s}(t)
$
be the size of subpopulation $\theta$ (independent of $t$ by
hypothesis). We do not make any specific hypothesis on the
ratios $\frac{\bar{M}_{\theta}^N}{N}$ as $N$ gets large (it may
be constant or not, it may tend to $0$ or not).

$\bullet$  \emph{Strategies and local interaction: } At time
slot $t,$ an ordered list $\calKK^N (t)$, of players in $ \{1,
2,\ldots,N\}$, without repetition, is selected randomly as
follows. First we draw a random number of players $K(t)$ such
that
 \ben
 \P(K(t)=k|M^N(t)=\vec{m})=J^N_{k}(\vec{m})
 \een
where the distribution $J^N_{k}(\vec{m})$ is given for any $N$,
$\vec{m} \in \{0,\frac{1}{N},\frac{2}{N},\ldots,1\}^{\sharp
\calX}$. Second, we set $\calKK^N$ to an ordered list of $K(t)$
players drawn uniformly at random among the
 $
 N (N-1) ... (N-K(t)+1)
 $
 possible ones. By abuse of notation we write $j\in \calKK^N(t)$ with the meaning
 that $j$ appears in the list $\calKK^N(t)$.

Each player such that $j\in \calKK^N (t)$ takes part in a
one-shot event at time $t$, as follows. First, the player
chooses an action $a$ in the finite set $\calA$ with
probability $u_{\theta}(a| s)$ where $(\theta, s)$ is the
current player state. The stochastic array $u$ is the strategy
profile of the population, and $u_{\theta}$ is the strategy of
subpopulation $\theta.$ A vector of probability distributions $u$ which depend only on the type of the player and its internal state is called {\it stationary strategy}.

Second, say that $\calKK^N(t)=(j_1, \ldots, j_k)$. Given the
actions $a_{j_1}, ...,a_{j_k}$ drawn by the $k$ players, we
draw a new set of internal states $(s'_{j_1}, ...,s'_{j_k})$
with probability
$ L^N_{\underline{\theta};\underline{s};\underline{a};\underline{s}'}
(k,\vec{m}) $,
 \bearn
\mbox{where }
\underline{\theta} = ( \theta_{j_1}, ...,\theta_{j_k}),\
\underline{s}=(s_{j_1}, ...,s_{j_k})
 \\
 \underline{a}=(a_{j_1}, ...,a_{j_k}),\
 \underline{s}'=(s'_{j_1},...,s'_{j_k})
 \eearn
 Then the collection of $k$ players makes one synchronized
transition, such that
 \ben
 S^N_{j_i}(t+\frac{1}{N})=s'_{j_i} \;\; \; i=1,\ldots,k
 \een
Note that $S^N_{j}(t+\frac{1}{N})=S^N_j(t)$ if $j$ is not in
$\calKK^N(t)$.

It can easily be shown that this form of interaction has
following properties: (1) $X^N$ is Markov and (2) players can
be observed only through their state.

The model is entirely specified by the probability
distributions $J^N$, the Markov transition kernels $L^N$ and
the strategy profile $u$. In this paper, we assume that $J^N$
and $L^N$ are fixed for all $N$, but $u$ can be changed and
does not depend on $N$ (though it would be trivial to extend
our results to strategies that depend on $N$, but this appears
to be unnecessary complication). We are interested in large
$N$.

It follows from our assumptions that
 \begin{enumerate}
 \item $M^N(t)$ is Markov.
 \item for any fixed $j\in\{1,\ldots,N\}$, $(X^N_j(t),
     M^N(t))$ is Markov. This means that the evolution of
     one specific player $X^N_{j} (t)$ depends on the other
     players only through the occupancy measure $M^{N}(t).$
 \end{enumerate}
\subsection{Payoffs}
We consider two types of instantaneous payoff and one
discounted payoff:

%\begin{itemize}
$\bullet$ %\item
{\it Instant Gain:} This is the random gain
    $G^N_j(t)$ obtained by one player whenever it is
    involved in an event at time $t$. We assume that it depends on
    this player's state just before the event and just after the event, the
    chosen action, and on the states and actions of all
    players involved in this event.
Formally, if player $j\in \calKK^N(t)$
 \ben G^N_j(t)=g^N(x_j,
a_j,x'_j,x_{\calKK^N(t)\backslash j}, a_{\calKK^N(t)\backslash
j},x'_{\calKK^N(t)\backslash j})
 \een
where $x_j=X^N_j(t)$, $a_j$ is the action chosen by player $j$,
$x'_j=X^N_j(t+\frac{1}{N})$, $x_{\calKK^N(t)\backslash j}$
[resp. $x'_{\calKK^N(t)\backslash j}$] is the list of states at
time $t$ [resp. at time $t+\frac{1}{N}$] of players other than
$j$ involved in the event, $a_{\calKK^N(t)\backslash j}$ is the
list of their actions and $g()$ is some non random function
defined on the set of appropriate lists. Whenever $j$ is not in
$\calKK^N(t)$, $G^N_j(t)=0.$ We assume that $G^N_j(t)$ is
bounded, i.e. there is a non random number $C_0$ such that,
with probability 1: for all $j, t$: $|G^N_j(t)|\leq C_0$

$\bullet$ %\item
{\it Expected Instant Payoff:} It is defined as the
    expected instant gain of player $j$, given the state $x$ of $j$
    and the population profile $\vec{m}$. By
    our indistinguishability assumption, it does not depend
    on the identity of a player, so we can write it as
    \bearn r^N(u,x,\vec{m}) :=\E\left(G^N_j(t)\right.\left|X^N_j(t)=x,M^N(t)=\vec{m}
    \right) \eearn
Note that this conditional expectation contains the case when
$j$ is not in $\calKK^N(t)$, i.e. when  $G^N_j(t)= 0$.

$\bullet$ %\item
{\it Discounted  Long-Term Payoff:} It is defined as the
expected discounted long term payoff of one player, given the
initial state of this player and the population:
$ \bar{r}^N(u;x,\vec{m}) :=  $
$$
\E ( \sum_{t=0 \mbox{ step }
1/N}^{\infty} e^{-\beta t}G_j^N(t)
  |
 X_j(0)=x,M^N(0)=\vec{m}
  )
$$
 where $\beta$
is a positive parameter (existence follows from the boundedness
of $G^N_j$). The fact that it does not depend on the identity
$j$ of the player, but only on its initial state $x$ and the
initial population profile $\vec{m}$, follows from the
indistinguishability assumption.

We defined the  Discounted  Long-Term Payoff in terms of the
instant gain, as this is the most natural definition. The
following proposition shows that the alternative definition, by
means of the expected instant payoff, is equivalent.
 \begin{subprop} \label{prott1}
For all player state $x$ and population profile $\vec{m}$
 \bearn
  \bar{r}^N(u;x,\vec{m}) &=& \E( \sum_{t=0 \mbox{ step }
1/N}^{\infty} e^{-\beta t}r^N(u,X^N_j(t),\vec{M}^N(t))
 \\
 &&|
 X_j(0)=x,M^N(0)=\vec{m}
 )
 \eearn
 \end{subprop}

\subsection{Focus
on One Single Player}
\label{sec-coi}

We are interested in the following special case (here we make
the dependency on the strategy explicit). There are two types
of players, i.e. $\Theta= 2$. There is exactly one player (the
player of interest) with type $1$. All other players have type
2. In this case we use the notation
 $R^N(u_1,u_2; s, \vec{m})$
 for the discounted long-term payoff obtained by the player in
 type $0$, when her strategy is $u_1$ and all other players's
 strategy is $u_2$, given that this player's initial internal
 state is $s$ and the initial type 2 subpopulation profile is $\vec{m}$.
Note that \ben
 R^N(u_1,u_2; s, \vec{m}) =  \bar{r}^N(u_1,u_2; (1,s), \vec{m}')
 \een
 with $m'_{1,s'}=\frac{1}{N}\ind{s=s'}$ and $m'_{2,s'}=m_{2,s'}$
 for all $s' \in \calS$.

\subsection*{ Mean Field Markov Game}
Player $j$ may choose a strategy $u_j$ which laws depends on its type and its own-internal state.
We look for a (Nash) equilibrium $u$ such that if
all players use $u$ then no player has an incentive to deviate
from $u$. For any finite $N$ one can map this into a standard
Markov game. This is true for both the case where the
number of players is known and in the case it is unknown
when taking a decision. Therefore we know that a stationary
equilibrium exists in the discounted case. A stationary equilibrium is solution of the fixed point equation: \[ \forall j,\
u_{j,\theta} \in \arg\max_{v_{j,\theta}} R^N ( v_{j,\theta} , u_{-j} ; s , m )
\]
 By assuming indistinguishability per type we can show that
a stationary equilibrium exists which is a solution
of the fixed point equation
\[ \forall \theta,
u_{\theta} \in \arg\max_{v_{\theta}} R^N (v_{\theta} , u ; s , m )
\]

Note that the mean field optimality here refers to the maximization of $R^n(u,u,s,m)$ over symmetric and stationary strategies. It is not necessarily optimal in the global sense.

\subsection*{ Mean Field Markov  Team} We wish to
find a stationary $u$ that maximizes $R^N$ averaged over
all players. \[
u=(u_{1},\ldots,u_{\Theta}) \in \arg\max_{v} R^N ( v ; s , m )
\]

\section{Main Results} \label{secmain}
\subsection{Scaling Assumptions}
\label{sec-scaling} We are interested in the large $N$ regime and
obtain that, for any fixed $j$, $(X_j^N,M^N)$ converges weakly to a
simple process. This requires the weak convergence of $M^N(0)$ to
some $\vec{m}_0$.

We assume that the parameters of the model and the payoff per
time unit converge as $N \to {\infty}$, i.e.
 \be
 \bracket{
 J^N_{k}(\vec{m})  \to   J_{k}(\vec{m})
\\
  L^N_{\underline{\theta};\underline{s};\underline{a};\underline{s'}}(k,\vec{m})
   \to
     \;\;\;\;\;\; L_{\underline{\theta};\underline{s};\underline{a};\underline{s'}}(k,\vec{m})
  \\
  r^N(u,x,\vec{m})  \to r(u,x,\vec{m})
 }
  \label{eq-scaling}
 \ee
Our main scaling assumption is
\begin{description}
\item[H1] $\sum_{k} k^2 J_{k}(\vec{m}) < \infty$ for all
    $\vec{m}\in \Delta$. This ensures that the second
    moment of the number of players involved in an event
    per time slot is bounded.
\end{description}
Note that H1 excludes the case where the number of players
involved in an event per time slot scales like $N$ (i.e.
synchronous transitions of all players at the same time). There
may be large $N$ asymptotic results for such cases
\cite{boudec2007gmf} but the limit is not given by an ODE. In
contrast, H1 is automatically true if the number of players
involved in an event per time slot is upper bounded by a non
random constant.
We also need some technical assumptions, which are usually true
and can be verified by inspection.
\begin{description}
\item[H2] $\sum_{k}J_k(\vec{m})
            > 0$ for all $\vec{m}\in \Delta$ ($\Delta$ is the simplex
        $\{\vec{m}: m_{\theta, s} \geq 0, \sum_{\theta, s}
        m_{\theta, s}=1 \}$). This
            ensures that the mean number of players
            involved in an event per time slot,
            $\sum_{k\geq 0} k J_k(\vec{m})$ is non zero.
\end{description}

 Define the \emph{drift} of $M^N(t)$ as \ben
        \vec{{f}}^{N}(u,\vec{m}) =
        \mathbb{E}\left({M}^{N}(t+\frac{1}{N})-{M}^{N}(t) |
        {M}^{N}(t)=\vec{m}\right) \een
Note that we make explicit the dependency on the strategy $u$
but not on $J$ and $L$, assumed to be fixed.

It follows from
our hypotheses that
  \be \lim_{N\to {\infty}}N f^N(u,\vec{m}) := f(u,\vec{m})
\label{eq-cv-drift} \ee exists.

   \begin{description}
    \item[H3]We assume that the convergence in
        \eref{eq-cv-drift} is uniform in $\vec{m}$ and the
        limit is Lipschitz-continuous in $\vec{m}$. This is
        in particular true if one can write, for every
        strategy $u$,
        $f^N(u,\vec{m})=\frac{1}{N}\phi_u(\frac{1}{N},\vec{m})$,
        with $\phi_u$ defined on $[0,\epsilon]\times
        \Delta$ where $\epsilon >0$ and $\Phi_u$ is
        continuously differentiable.
        \item[H4] $\P(X^N_j(t+1/N)=y|X^N_j(t)=x,M^N(t)=m,
            M^N(t+1/N)=m')$  converges uniformly in
            $\vec{m},\vec{m}'$ and the limit is
            Lipschitz-continuous in $\vec{m},\vec{m}'.$
            This is in particular true if one can write,
            for every strategy $u,$    as
            $\xi_{u,x;y}(1/N,m,m').$ with $\xi$ defined on
            $[0,1]\times \Delta\times  \Delta $ and
            $\xi_{u,x;y}$ is continuously differentiable.
\end{description}

Our model satisfies the assumptions in \cite{JYLB08}, therefore we
have the following result:
\begin{subthm}[\cite{JYLB08}] \label{Gra2} Assume that
$\lim_{N\longrightarrow\infty}M^N(0)=\vec{m}_0$ in probability.
For any stationary strategy $u,$ and any time $t,$ the random
process $M^N(t)=\frac{1}{N}\sum_{j=1}^N\delta_{X_j^N(t)}$
converges in distribution to the (non-random) solution of the
ODE
 \be
 \dot{\vec{m}}(t)=f(u,\vec{m}(t))
 \label{eq-ode}
 \ee
 with initial condition $\vec{m}_0$.
\end{subthm}
%In the case with $\Theta=2$ types, we call
%${R}(u_1,u_2;\alpha,\vec{m})$ the discounted long term
%payoff.

\subsection{Convergence results}
We focus on one player, without loss of generality we can call
her player $1$, and consider the process $(X^N_1, M^N)$. For
any finite $N$, $X^N_1$ and $M^N$ are not independent, however
in the limit we have the following:
\begin{subthm}
\label{prop2} Assume that
$\lim_{N\longrightarrow\infty}M^N(0)=\vec{m}_0$ and
$\lim_{N\longrightarrow\infty}X^N_1(0)=x_0=(\theta_1,s_0)$
 in
probability. The discrete time process $(X^N_1(t), M^N(t))$
defined for $t\in \frac{\Nats}{N}$, converges weakly to the
continuous time jump and drift process $(X_1(t), \vec{m}(t))$,
where $\vec{m}(t)$ is solution of the ODE \eref{eq-ode} with
initial condition $\vec{m}_0$ and $X_1(t)$ is a continuous
time, non homogeneous jump process, with initial state $x_0$.
The rate of transition of $X_1(t)$ from state
$x_1=(\theta_1,s_1)$ to state $x_1'=(\theta_1, s_1')$ is
 \ben A(x_1,x_1';\vec{m}(t),u)=\sum_{k\geq 1}J_k(\vec{m})A_k(s_1,s_1';\vec{m}(t),u)\een
with $ A_k(s_1,s_1';\vec{m}(t),u)= $
\[
\sum_{\underline{\theta};\underline{s};\underline{a};\underline{s}'}
 L_{\theta_1,\underline{\theta};s,\underline{s};\underline{a};s',\underline{s}'}(k,\vec{m}(t))
 \prod_{j=1}^k u_{\theta_j}(a_j|s_j)
 \prod_{j=2}^k m_{\theta_j,s_j}(t)
 \]
\bearn \mbox{where }\underline{\theta} = ( \theta_{2}, ...,\theta_{k}),
\underline{s}= (s_{2}, ...,s_{k})\\
 \underline{a}=(a_{1}, ...,a_{k}),
 \underline{s}'=(s'_{2},...,s'_{k})
 \eearn
 \end{subthm}

Note that, contrary to results based on propagation of chaos,
we do not assume that the distribution of player states at time
$0$ is exchangeable. In contrast, we will use
Theorem~\ref{prop2} precisely in the case where player $1$ is
different from other players. Theorem~\ref{prop2} motivates the
following definition.

\begin{subdef} To a game as defined in \sref{sec-model} we
associate a ``Macroscopic Mean Field Markov Game",
defined as follows. There is one player, (player 1), with state
$X_1(t)$ and a population profile $\vec{m}(t)$. The initial
condition of the game is $X_1(0)=x$, $\vec{m}(0)=\vec{m}_0$.
The population profile is solution to the ODE~(\ref{eq-ode})
and $X_1(t)$ evolves as a jump process 
Theorem~\ref{prop2}.

Further, let $\bar{r}(u;x,\vec{m})$ be the discounted long-term
payoff of player 1 in this game, given that $X_1(0)=x$ and
$\vec{m}(0)=\vec{m}_0$, i.e. $\bar{r}(u;x, \vec{m})=$
 \ben
\E\left(\int_0^{\infty}e^{-\beta t}
 r(u,X_1(t),m(t))|X_1(0)=x,\vec{m}(0)=\vec{m}_0\right)
\een

We also consider, as in \sref{sec-coi}, the case with
$\Theta=2$ types and define by analogy $R(u_1,u_2; s, \vec{m})$
as the discounted long-term payoff when player $1$ starts in
state $s$ and the population profile starts in state $\vec{m}$,
with player 1 using strategy $u_1$ and other players strategy
$u_2$.
\end{subdef}

In order to exploit the convergence in distribution of the process
focused on one player, we need that the payoff be continuous in the
topology of this convergence. This is stated in the next theorem.
\begin{subthm}
Let $E=\calS\times\Delta$ and $D_E[0,{\infty})$ the set of cadlag
functions from $[0,{\infty})$ to $\Reals$, equipped with Skorohod's
topology. The mapping
  \bearn
  D_E[0,{\infty}) &\to & \Reals\\
  (s,m)&\mapsto&
  \int_0^{\infty}e^{-\beta t}r(u,s(t),m(t))\ dt
  \eearn
is continuous. \label{theo-rew-cont}
\end{subthm}

Using Theorem~\ref{prop2} and Theorem~\ref{theo-rew-cont} we
obtain the following, which is the main result of this paper.
It says that when $N$ goes to infinity, the Mean Field Markov Game with $N(t)$ of players becomes equivalent to
the associated Macroscopic Mean Field Markov Game.
This reduces any multi-player problem into an effective
one-player problem facing an evolving aggregative object.

\begin{subthm}[Asymptotically equivalent game] \label{mainthm}
When $N$ goes to infinity we have (a) the discrete time process
$X^N_1$  converges in
      distribution to the continuous time process $X_1$
  (b) $\bar{r}^N(u;x,\vec{m})\to \bar{r}(u;x, \vec{m})$
  and (c)  $R^N(u_1,u_2; s, \vec{m})\to R(u_1,u_2; s,
      \vec{m})$
 \end{subthm}

\subsection{Case with Global Attractor}
\label{sec-glo} Assume that, for some strategy $u$, the ODE
(\ref{eq-ode}) has a global attractor $\vec{m}^*$  (this may or may
not hold, depending on the ODE). If in addition the model with $N$
players is irreducible, with stationary probability distribution
$\varpi^N$ for $M^N$, then $\lim_{N\longrightarrow \infty}\
\varpi^N= \delta_{\vec{m}^*}$ where $\delta_{\vec{m}^*}$ is the
Dirac mass at $\vec{m}^*$ (follows from \cite{JYLB08}). i.e. the
large time distribution of $M^N(t)$ converges, as $N\to \infty$, to
the attractor $\vec{m}^*$.

Also,  $(X^N_j(t), M^N(t))$ converges to a continuous time,
homogeneous Markov jump process with time-independent
transition matrix:
\ben
 A(x_1, x'_1; u) = \sum_{k\geq 1}J_k(\vec{m})A_k(s_1,s_1';\vec{m}^*,u)
 \een
Assume that the transition matrix $A(x_1, x'_1; u)$ is also
irreducible and let $\pi()$ be its unique stationary
probability. Also let $\pi^N$ be the first marginal of the
stationary probability of $(X^N_1,M^N)$. It is natural in this
case to replace the definition of the long term payoffs
$R^N(u_1,u_2;s,\vec{m})$ and $R^N(u_1,u_2;s,\vec{m})$ by their
stationary counterparts
 \bearn
 R_{st}^N(u_1,u_2)&:=&\sum_s \pi^N(s)R^N(u_1,u_2;s,\vec{m}^*)
 \\
R_{st}(u_1,u_2)&:=&\sum_s \pi(s)R(u_1,u_2;s,\vec{m}^*)
 \eearn

\subsection{Single player per type selected per time slot}
Consider the special case where at each time slot, only one
player per type between the $N$ is randomly selected and has a
chance to change its action, i.e. $\sharp\calKK^N=1$ w.p 1.

Thus H1 and H2 are automatically satisfied. The resulting ODE (see
\cite{BW03}) becomes

$$ \frac{d}{dt}m_{x}(t)=\sum_{x'}m_{x'}L_{x',x}(\vec{m}, u,\Theta)-m_{x}\sum_{x'}L_{x, x'}(\vec{m}, u,\Theta)$$

The term $\sum_{x'}m_{x'}L_{x',x}(\vec{m}, u,\Theta)$ is the {\it
incoming flow} in to $x$ and the {\it outgoing flow} from $x$ is
$m_{x}\sum_{x'}L_{x, x'}(\vec{m}, u,\Theta).$

We then obtain a large class of  {\it state-dependent} evolutionary
game dynamics. Note that in general the trajectories of the mean
dynamics need not to converge. In the case of single player selected
in each time slot of $1/N$ and linear transition in $m,$ the time
averages under the replicator dynamics converge its interior rest points or
the boundaries of the simplex.

\subsection{Equilibrium and optimality}
Let $\calU_s$ be the set of strategies. Consider the optimal
control problems
$$ (OPT_N)\  \left\{\begin{array}{c}\mbox{ Maximize}\ R^N (u,u;s,\vec{m}_0)\\ \mbox{s.t}\ u\in \calU_s \end{array}\right.$$
$$ (OPT_{\infty})\  \left\{\begin{array}{c}\mbox{ Maximize}\ R (u,u;s,\vec{m}_0)\\ \mbox{s.t}\ u\in \calU_s \end{array}\right.$$
The strategy $u$ is an $\epsilon-$optimal strategy for the
$N$-optimal control problem if $$\ R^N(u,u;s,\vec{m}_0)\geq
-\epsilon+\sup_{v}\ R^N(v,v;s,\vec{m}_0).$$

Also consider the fixed-point problems
$$ (FIX_N)\ \left\{\begin{array}{c}\mbox{find}\ u\in \calU_s \mbox{ such that}\\ u\in\arg\max_{v\in\calU_s}\{R^N(v,u;s,\vec{m}_0)\}  \end{array}\right.$$
$$ (FIX_{\infty})\ \left\{\begin{array}{c}\mbox{find}\ u\in \calU_s \mbox{ such that}\\ u\in\arg\max_{v\in\calU_s}\{R(v,u;s,\vec{m}_0)\}  \end{array}\right.$$
A solution to $(FIX_N)$ or $(FIX_{\infty})$ is  a ( Nash)
equilibrium. We say that $u$ is an $\epsilon-$equilibrium for the
game with $N$ [resp. $N \to \infty$] players  if
$R^N(u,u;s,\vec{m}_0)\geq \sup_{v}
R^N(v,u;s,\vec{m}_0)-\epsilon$ [resp. $R(u,u;s,\vec{m}_0)\geq
\sup_{v} R(v,u;s,\vec{m}_0)-\epsilon$].

Note that the definition of equilibrium and optimal strategy
may depend on the initial conditions. If, for any $u \in
\calU_s$, the hypotheses in \sref{sec-glo} hold, then we may
relax this dependency.

\begin{subthm}[Finite $N$] \label{prop3}
For every discount factor $\beta>0$ the optimal control problem $(OPT_N)$ (resp. the fixed-point problem $(FIX_N)$) has at least one $0-$optimal strategy (resp. $0-$equilibrium). In particular, there a
  $\epsilon_N$-optimal strategy (resp. $\epsilon_N-$equilibrium)  with $\epsilon_N\longrightarrow 0$.
\end{subthm}

\begin{subthm}[Infinite $N$] \label{propt4}
 Optimal strategies (resp. equilibrium strategies) exist in the limiting regime when $N\to \infty$  under uniform convergence and continuity of $R^N\to R.$ Moreover, if $\{U^N\}$ is a sequence of $\epsilon_N-$optimal strategies (resp. $\epsilon_N-$equilibrium strategies) in the finite regime with $\epsilon_N\longrightarrow \epsilon$, then, any limit of subsequence $U^{\phi(N)}\longrightarrow U$ is an $\epsilon-$ optimal strategy (resp. $\epsilon-$equilibrium) for game with infinite $N.$
\end{subthm}

\subsection{Mean field equilibrium}\label{sec37}

Each generic player $1$ with strategy $v_1$ optimizes its own long-term payoff under the behavior of its own-internal state which is a continuous time Markov jump process driven by $A(x_1,x_1';v_1,\vec{m}(t),u)$ and the behavior of $\vec{m}(t)$ is given by the controlled ODE under the strategy $u.$

It is important to notice that at the infinite population limit the mean field limit dynamics does not depend on $v_1.$ This can be easily seen from the fact that the effect of a single player is in order of $\frac{1}{N}.$ When $N\longrightarrow +\infty,$ the effect becomes negligible with the respect to the mass. However, $v_1$ can be a big effect in the rate transition of that player via  $A(x_1,x_1';v_1,\vec{m}(t),u).$

The consistency between the individual state transition and the fraction of players per state needs to be checked.

We say that the pair $(u^*_t,\vec{m}^*(t))$ is  a mean field equilibrium if
$\{u^*_t\}_{t\geq 0}$ is a mean field response to the individual dynamic optimization where $\vec{m}^*(t)$ is the mean field  at time $t$ and $u^*_t$ produces the mean field  i.e $\vec{m}[u^*,\vec{m}_0](t)=\vec{m}^*(t).$

If $(v^*,u^*,m^*)$ satisfies the following equation
$$\left\{
\begin{array}{c}
\!\!\!\beta v_{\theta,t}(s,m)=\displaystyle{\sup_{u_{\theta}}}\left\{ r_{\theta}(y_{\theta},u_{\theta},\vec{m}(t)){+}\sum_{s'}A(s,s';\vec{m}(t),u) v_{\theta,t}(s'_{\theta},\vec{m}(t))\right\}+f(u^*,m).\partial_mv_{\theta,t}
\\
m_{\theta}(t)=m_{\theta,0}+\int_0^t {f}_{\theta}(u^*_{t'},\vec{m}(t'))\ dt'\\
\ m(0)=m_0\in\ \Delta(\mathcal{X}),\
\theta\in\Theta.
\end{array}
\right.$$
and the strategy $$u^*_{\theta, t}\in\arg\max\left\{ r_{\theta}(y_{\theta},u_{\theta},\vec{m}(t)){+}\sum_{s'}A(s,s';\vec{m}(t),u) v_{\theta,t}(s'_{\theta},\vec{m}(t))\right\},$$
then, one gets a mean field equilibrium.
The problem becomes
$\max_{v_1\in\calU_s}\{R(v_1,u;x_1,\vec{m}_0)\} $ subject to  the transitions $A(x_1,x_1';v_1,\vec{m}(t),u)$ and the ODE.
\section{Illustrating example} \label{hd}

We present in this section an example of a dynamic version of the
Hawk and Dove problem where each individual has three energy levels.
We derive the mean field limit for the case where all users
follow a given policy and where possibly one player deviates.
We then further simplify the model to only two energy states
per player. In that case we are able to fully identify
and compute the equilibrium in the limiting Mean Field Markov Game.
Interestingly, we show that the ODE converges to a fixed point
which depends on the initial condition and the policy.

Consider an homogenous population of $N$ animals. An
animal plays the role of  a  {\it player}. Occasionally two
animals find themselves in competition on the same piece of
food.
Each animal has three states $x=0,1,2$
which represents its energy level. An animal can adopt
an aggressive behavior (Hawk) or a
peaceful one (Dove, passive attitude). At the state $x=0$ there is no action.
We describe the fitness
of an animal (some arbitrary player) associated with the possible
outcomes of the meeting as a function of the decisions taken by
each one of the two animals. The fitnesses represent the following:
\begin{itemize}
\item  An encounter Hawk-Dove or Dove-Hawk results in zero fitness to the
Dove and in $\bar{v}$ of value for the Hawk that gets all the food
without fight. The state of the Hawk (the winner) is incremented $a=\ind{\{x'_H=\min(x_H+1,2)\}}$ and the state of the Dove is $b=\ind{\{x'_D=\max(x_D-1,0)\}}.$

\item
An encounter Dove-Dove results in a peaceful, equal-sharing of the
food which translates to a fitness of $\frac{\bar{v}}{2}$ to each
animal and the state of each animal change with the sum of the two
distributions $ \frac{1}{2}a+\frac{1}{2}b$
\item
An encounter Hawk-Hawk results in a fight in which with $p=1/2$
chances, one (resp. the other) animal obtains the food but also
in which there is a positive probability for each one of the
animals to be wounded $1/2$ . Then the fitness of the animal 1 is
$\frac{1}{2}(\bar{v}-c)+\frac{1}{2} (-c)=\frac{1}{2}\bar{v}-c,$ where the
$-c$ term represents the expected loss of fitness due to being
injured.
\end{itemize}
{
$$
 \begin{array}{|ccc|}
\hline
i \backslash j &  (g^N_i,g^N_j)  & X^{N}_i(t+\frac{1}{N}),X^{N}_j(t+\frac{1}{N})
\\ \hline
D - D & (\frac{\bar{v}}{2},\frac{\bar{v}}{2}) & \frac{1}{2} \delta_{\min(x_1-1,0),\max(x_2+1,2)}\\ & & +\frac{1}{2} \delta_{\max(x_1+1,2),\min(x_2-1,0)}                                                                                \\ \hline
D - H & (0, v) &  (\min(x_1-1,0),\max(x_2+1,2))\\ \hline
H - H & \frac{1}{2}v-c &  \frac{1}{2} \delta_{\min(x_1-1,0),\max(x_2+1,2)}\\ & & +\frac{1}{2} \delta_{\max(x_1+1,2),\min(x_2-1,0)}\\ \hline
\end{array}
$$
} The vector of frequencies of states at time $t$   is given by
$M^{N}_{x}(t)=\frac{1}{N}\sum_{j=1}^N\ind{\{X^N_{j}(t)=x\}}$
for $x=0,1,2$ and the action set is $A_x=\{H,D\}$ in each state
$x\neq 0,$ $A_0=\{\}.$

The assumptions in Section \ref{secmain} are satisfied (pairwise interaction, $\sharp \calKK^N(t)={2}$) and the occupancy measure  $M^{N}(t)$ converges to $m(t).$

\subsection{ODE and Stationary  strategies}
Consider the following fixed parameters
$\mu_1=L_{0,1},\ \mu_2=L_{0,2}.$ The population profile is
denoted by $\vec{m}=(m_{0},m_{1},m_{2})$ and the stationary strategy
is described by the parameters $v_1,v_2$
where $v_1:=u(H|1),\ v_2=u(H|2)$ \\
{
$$
\begin{array}{l}
\dot{m}_{2}=  m_0 L_{0,2}+m_1 L_{1,2}(u,m) -m_2 L_{2,1}(u,m)
\\
\dot{m}_{1}=  m_0 L_{0,1}+m_2 L_{2,1}(u,m)\
-m_1L_{1,2}(u,m)-m_1 L_{1,0}(u,m))\\
\dot{m}_{0}= m_1 L_{10}(u,m)-(\mu_1+\mu_2)m_0
\end{array}
$$
}
\\
where
{
$ {L_{12}(u,m)= } $
\begin{eqnarray*}
&& m_0+v_1\left(\frac{v_1 m_1}{2}+(1-v_1)m_1+\frac{v_2 m_2}{2}+(1-v_2)m_2\right)
\\
& &+ (1-v_1)\left( \frac{(1-v_1)m_1}{2}+\frac{(1-v_2)m_2}{2}\right)
\\
\lefteqn{ L_{2,1}(u,m) = v_2\left(  \frac{v_1 m_1}{2}+\frac{v_2 m_2}{2}   \right) } \\
 &  &+ (1-v_2)\left( \frac{(1-v_1) m_1}{2}+ v_2 m_2+\frac{(1-v_2)m_2}{2}\right)
\\
\lefteqn{ L_{10}(u,m)  :=  v_1\left(\frac{v_1 m_1}{2}+\frac{v_2 m_2}{2} \right) }
\\
&+& (1-v_1) \left(v_1 m_1+\frac{(1-v_1) m_1}{2}+v_2 m_2+\frac{(1-v_2) m_2}{2}
\right)  ,
\end{eqnarray*}
}

%Our  numerical experiment studies the behavior and illustrates the evolution of frequencies of states.
%
%{\bf to be added}
For $\calKK^N=\{ j_1,j_2\},$
$x'_j,x_i\in\{0,1,2\},$
\begin{eqnarray}
\frac{d}{dt}m_{x} &=&  \sum_{x_1,x_2,x_2'}
 m_{x_1}m_{x_2}L_{x_1,x_2;x,x_2'}(u,\vec{m})
 \nonumber\\
 &&+
\sum_{x_1,x_2,x_1'}
 m_{x_1}m_{x_2}L_{x_1,x_2;x_1',x}(u,\vec{m})
\nonumber\\
&&- m_x \sum_{x_2,x_1',x_2'}
 m_{x_2}L_{x,x_2;x_1',x_2'}(u,\vec{m})
 \nonumber\\ & &- m_x \sum_{x_1,x_1',x_2'}
 m_{x_1}L_{x_1,x;x_1',x_2'}(u,\vec{m}) \nonumber
\end{eqnarray}

\subsection{Computation of $R(u_1,u_2; s,\vec{m})$.}
 We want to compute the value $$V(u_1,u_2,x,m):= \mathbb{E}_{x}\int_{0}^{\infty}e^{-\beta t}r(u_1,u_2,x(t),m(t))\ dt$$
$$s.t. \ \dot{m}(t)=f(u_2,m(t)), m(0)=m_0, \ x(0)=x.$$
\begin{eqnarray}
V(u_1,u_2,x,m)&=& \mathbb{E}_{x}\int_{0}^{\Delta}e^{-\beta t}r(u_1,u_2,x(t),m(t))\ dt \nonumber\\
& & +\mathbb{E}_{x}\int_{\Delta}^{\infty}e^{-\beta t}r(u_1,u_2,x(t),m(t))\ dt \nonumber\\
%&=&  \mathbb{E}_{x}\int_{0}^{\Delta}e^{-\beta t}r(u_1,u_2,x(t),m(t))\ dt \nonumber\\
%\nonumber &+& \mathbb{E}_{x}e^{-\beta \Delta}\int_{\Delta}^{\infty}e^{-\beta (t-\Delta)}r(u_1,u_2,x(t),m(t))\ dt\\ \nonumber &= &
% \mathbb{E}_{x}\int_{0}^{\Delta}e^{-\beta t}r(u_1,u_2,x(t),m(t))\ dt \nonumber\\ \nonumber &+& \mathbb{E}_{x}e^{-\beta\Delta} \mathbb{E}_{x(\Delta)}\int_{0}^{\infty}e^{-\beta t}r(u_1,u_2,x(t),m(t))\ dt\\
\nonumber & =&
 \mathbb{E}_{x}\int_{0}^{\Delta}e^{-\beta t}r(u_1,u_2,x(t),m(t))\ dt  \\ \nonumber  & & +\mathbb{E}_{x}e^{-\beta\Delta} V(u_1,u_2,x(\Delta),m(\Delta))
\end{eqnarray}
This implies that
\begin{eqnarray} \nonumber
%0 %&=& \frac{1}{\Delta} \mathbb{E}_{x}\left\{ \int_{0}^{\Delta}e^{-\beta t}r(u_1,u_2,x(t),m(t))\ dt \right.\\ \nonumber & + & \left.
%e^{-\beta\Delta} V(x(\Delta))-V(x) \right\}\\ \nonumber
0& =& \mathbb{E}_{x}\frac{1}{\Delta} \int_{0}^{\Delta}e^{-\beta t}r(u_1,u_2,x(t),m(t))\ dt  \\ \nonumber & &+
\frac{e^{-\beta\Delta}-1}{\Delta} \mathbb{E}_{x} V(u_1,u_2,x(\Delta),m(\Delta))\\ & & + \frac{\mathbb{E}_{x} V(u_1,u_2,x(\Delta),m(\Delta))-V(u_1,u_2,x,m)}{\Delta}
\end{eqnarray}
Using Ito's formula and Lebesgue integration properties, we obtain that:
$\frac{\mathbb{E}_{x} V(u_1,u_2,x(\Delta))-V(u_1,u_2,x)}{\Delta}$ goes to $\sum_{x'}{D}_{m_{x'}}V(u_1,u_2,x')\frac{d}{dt}m_{x'}+jumps,$ where ${D}_{m_{x'}}V$ is the derivative of $V$ in a weak sense,   $\frac{e^{-\beta \Delta}-1}{\Delta}\longrightarrow -\beta ,$ and the term $$\mathbb{E}_{x}\frac{1}{\Delta} \int_{0}^{\Delta}e^{-\beta t}r(u_1,u_2,x(t),m(t))\ dt \longrightarrow r(u_1,u_2,x,m_0)$$ when $\Delta$ goes to zero, and the jump term is due to the changes in the process $x.$ The jump term is explicitly determined by the transitions rates which {\it contains} $u_1$ and the value $V$ as given section \ref{sec37}.
Thus, we obtain
 { \begin{eqnarray} \label{hjb}
\beta V(u_1,u_2,x,m)= r(u_{1,x},u_{2,x},x,m)+\sum_{x'}(D_{m_{x'}}V(u_1,u_2,x',m))f_{x'}(u_2,m)+jumps
\end{eqnarray}} where $u_{i,x}=u_i(H|x).$

The  optimality is then given by the
Hamilton-Jacobi-Bellman equation obtained by maximizing
the right-hand side of the equation (\ref{hjb}) over the action set.
{ \begin{eqnarray} \nonumber
\beta \Psi(x,m)= \max_{u_{1,x},u_{2,x}} \{r(u_{1,x},u_{2,x},x,m)+\sum_{x'}(D_{m_{x'}}\Psi(u_1,u_2,x',m))f_{x'}(u_2,m)+jumps\}
\end{eqnarray}} and optimality conditions of the best response to $u_2$ is given by
{ \begin{eqnarray} \nonumber
\beta \Phi(u_2,x)= \max_{a\in \{H,D\}} \{r(a,u_{2,x},x,m)+\sum_{x'}(D_{m_{x'}}\Phi(u_2,x'))f_{x'}(u_2,m)+jumps\}
\end{eqnarray}}
Note that in the global optimization case (under symmetry per class strategies) we can drop the jump terms by computing the expected social welfare (which do not depend on $x$ but depends on $m$). Hence the equation reduces to a similar one as in
\cite{gamenets}. Now, if we consider the individual optimization problem, there is a jump and drift term in the generator as it is usual in hybrid systems. 
Theses equations are in general difficult to solve and the solutions are not necessarily regular (e.g. viscosity solutions). Numerical approaches based on multi-grid techniques of Hamilton-Jacobi-Bellman-Issacs equations can be found \cite{Kushner}.
\subsection{The case of two energy levels}
In order to derive closed form expressions for solutions of our ODE,
we consider  two states, i.e., each animal has
two states $x=1, 2$ which represents its energy levels. Thus, the
ODE  % (\ref{ODEexample})
can be expressed as follows:
\begin{eqnarray}
\dot{m}_{2}(t)&=&  (1-m_2(t)) L_{1,2}(u,m) -m_2(t) L_{2,1}(u,m)
\end{eqnarray}
which can be rewritten as
\begin{eqnarray}
\dot{m}_{2}(t)&=&  a_1+a_2 m_{2}(t)+a_3 (m_2(t))^2
\end{eqnarray} with $a_1=1,$ $a_2=\frac{u_2}{2}-2<0,$ $a_3=\frac{1-u_2}{2}>0.$

Let $m[u,m_0](t)$ be the solution of the
ODE given $u$ and a initial distribution $m(0)=m_0.$
We distinguish two cases:
\begin{itemize}
\item[Case 1]   $u_2=1$ (fully aggressive when it is possible): the ODE becomes $ \dot{m}_2(t)=1-\frac{3}{2}m_2(t)$ and the solution has the form
     \begin{eqnarray}\label{mean1} m_{2}[1,m_0](t)=\frac{2}{3}[1-c_1 e^{-\frac{3}{2}t}]
      \end{eqnarray} with $c_1=1-\frac{3}{2}m_0$ and $m_1[u,m_0](t)=1-m_2[u,m_0](t)$

\item[Case 2] $u_2\neq 1,$ (less aggressive in state 2)
\begin{eqnarray}\label{mean2} m_{2}[u,m_0](t)=\gamma_{-}(u)+\frac{\gamma_{+}(u)-\gamma_{-}(u)}{1-c_2 e^{(\gamma_{+}(u)-\gamma_{-}(u))a_2 t}}\end{eqnarray}
\begin{eqnarray*}
\mbox{where }
 c_2 &=& 1+\frac{\gamma_{+}(u)-\gamma_{-}(u)}{m_{2}(0)-\gamma_{-}(u)},
\\
\gamma_{-}(u) &=& \frac{2-u_2/2-
(2+u_2^2/4)^{\frac{1}{2}}}{1-u_2}<1,
\\
\gamma_{+}(u) &=& \frac{2-u_2/2+(2+u_2^2/4)^{\frac{1}{2}}}{1-u_2} >1
\end{eqnarray*}
\end{itemize}
Note that in both cases there is a unique strategy-dependent global attractor.
$$ \lim_{t\longrightarrow \infty} m_2[u,m_0](t) = \left\{\begin{array}{c} \gamma_{-}(u) \ \mbox{if} \ u_2\neq 1\\   2/3 \ \mbox{if}\ u_2=1
\end{array}\right.$$

The expected instant payoff of a player using the stationary strategy $v$  when the
population profile is $m[u,m_0](t)$, is given by
{ \[
r(v,u,2,m[u,m_0](t))= v[\bar{v}-c m_{2}u_2]+(1-v)r(v,u,1,m[u,m_0](t))
\]
\[
r(v,u,1,m[u,m_0](t))=\frac{1}{2}(1-m_2[u,m_0](t)u_2)\bar{v}
\]
}
where $m_2[u,m_0](t)$ is given by (\ref{mean1}) (resp.
(\ref{mean2})) for $u_2=1$ (resp. $u_2\not=1$). Now, we can compute
explicitly the best response against $u$ for a given initial $m_0.$
Let
{

\[
\beta_2(u,2,m_0,t)=r(H,u,2,m[u,m_0](t))-r(D,u,2,m[u,m_0](t)).
\]
}
The best response, $ \mbox{BR}(x,u,m[u,m_0](t))$,
against $u$ at $t$ is
\[
\mbox{BR}(x,u,m[u,m_0](t)) =
\left\{
\begin{array}{lr}
\mbox{play Hawk if } \beta_2(u,x,m_0,t)>0\\
\mbox{play Dove if } \beta_2(u,x,m_0,t)<0
\end{array}
\right.
\]
This implies that it is better to play Hawk for
$\frac{\bar v}{2c}>\frac{\gamma}{1+\gamma}$ where $\gamma=\max(2/3,m_0).$
Since the solution of the ODE is strictly monotone in time for each stationary strategy, there is at most one time for which $\beta_2$ is zero.
It is easy to see that if $\frac{\bar v}{2c}>\frac{2}{3}$ then the strategy
which to play Hawk in state 2 and Dove in state 1 is an equilibrium.

\begin{figure}[htp]
\centering
  %\includegraphics[width=8cm,height=5cm]{HDgamenets1.pdf}\\
  \caption{Global attractor for $u_2=1$}\label{fighd1}
\end{figure}
\begin{figure}[htp]
\centering
  %\includegraphics[width=8cm,height=5cm]{HDgamenets2.pdf}\\
  \caption{Global attractor for $u_2=0.2$}\label{fighd2}
\end{figure}
\section{Concluding remarks} \label{concl}
The goal of this paper has been to develop mean field asymptotic of interactions with large number of players
using stochastic games. Due to the curse of the size of the population, the applicability of
atomic stochastic games  has been severely limited. As an alternative, we proposed a method for mean field Markov games
 where players make decisions
only based on their own state and the  global system state. We have showed under mild assumptions
 convergence results, where asymptotics were taken in
the number of players. The population state profile satisfies a system of non-linear ordinary differential equations.
We have considered very simple class of strategies that are functions only of  player's own state
and the population profile. We applied to Hawk-Dove interaction with several energy level and formulated the ODEs. We show that  the best response depends on the initial conditions.

\section*{Appendix}

\subsection*{Sketch of proof of Proposition \ref{prott1}}
Let $\tau^N$ be the first time after $t=0$ that $X^N_j(t)$ hits in some given state.
We show that
\be \bar{r}^N = \frac{1}{N} \esp{ \sum_{s=0 \mbox{ step }
1/N}^{\tau^N} e^{-\beta t} r^N\left(X^N_j(s),M^N(s)\right) }
\label{eq-1t}\ee

Define for $t\in \mathbb{N}/N$:
 \ben
 Z_t^N = \sum_{s=0 \mbox{ step }
1/N}^{t} e^{-\beta s} \left(G^N(s)-
r^N\left(X^N_j(s),M^N(s)\right)\right)
 \een
 we have, for $0 \leq s \leq t$:
 {  \bearn
 Q&:=& \espc\left(Z_t^N-Z_s^N |\calF^N_s\right)\\ &=&
\hspace{-0.5cm}
\sum_{u'=0 \atop \mbox{ step } 1/N}^{t} e^{-\beta u'} \espc\left(
G^N(u')-
r^N\left(X^N_j(u'),M^N(u')\right)| \calF^N_s\right) \eearn }
which can be written as
{ \bearn
\sum_{u'=0  \atop \mbox{ step } 1/N}^{t}
\hspace{-0.5cm}
e^{-\beta u'} \espc\left(
         \espc\left(G^N(u')-
         r^N\left(X^N_j(u'),M^N(u')\right)
         |
        \calF^N_{u'}\right)
       |
       \calF^N_s\right)\\
= 0
 \eearn }
 thus $Z_t^N$ is an $\calF^N_t-$ martingale. Now
 $\tau^N$ is a stopping time with respect to the
 filtration $\calF^N_t$ thus, by Doob's stopping time theorem:
 $
 \esp{Z_{t \wedge \tau^N}^N}=\esp{Z_{0 \wedge \tau^N}^N}=0
 $
Further, $
Z_{t \wedge \tau^N}^N \leq K |\tau^N|$ for some
constant $K$. Since $\tau^N$ is almost surely finite
and has a finite expectation, we can apply dominated
convergence (with $t\to {\infty}$) and obtain
$ \esp{Z_{\tau^N}^N}=0.$

\subsection*{Sketch of Proof of Theorem \ref{prop2}}

To prove the weak convergence of $Z^N,$ we check the following steps:
Without loss of generality, we took the set of states as $\calS=\{0,1,2,\ldots,\sharp\mathcal{S} \}$
$X_j^N$ has a jump $r$ with probability $$q^N_{i,i+r}(M^N(k))=\frac{1}{N}L^N_{i,i+r}(M^N(k),u))$$ and $M^N$ is the continuous process with drift $f^N.$
\begin{itemize}
    \item We introduce of $\tilde X^N_j$ by scaling with step size $\frac{1}{N}.$ Then, $Z^N=(X^N,M^N)$ is approximate in some sense by a discrete time process $\tilde Z^N=(\tilde X^N,\tilde m^N)$ where $\tilde m^N(k)=m(\lfloor Nt\rfloor)$ $m$ solution of the ODE
        with $\tilde X_j^N$ is the discrete time jump process with transition matrix  $$q^N_{i,i+r}(\tilde m^N(k))=\frac{1}{N}L^N_{i,i+r}(m(\frac{k}{N}),u)).$$
We show that $ d(X^N_j,\tilde X^N_j)\longrightarrow 0$ for any compact of time intervals.
\item $$\tilde Z^N=(\tilde X^N,\tilde m^N)\Longrightarrow (\tilde X,\tilde m) $$
$M^N([Nt])\longrightarrow m(t).$
    We derive the weak convergence of $Z^N$ to $(X,m)$ where $m$ is deterministic and $X$ is random.

\end{itemize}

\noindent{ Approximation by a discrete time process}

The following lemma follows from the lemma 1 and 3 in Benaim and Weibull (2003,2008), in which we incorporate  behaviorial  strategies.
\begin{lem} \label{lemm}
For every $t>0$ there exists a constant $c$ such that  for every $\epsilon>0$ and $N$ large enough one has
{
$$P(\sup_{0\leq\tau\leq T} ||M^N(\tau)-m(\tau)|| >\epsilon | \ {M}^{N}(0)=m_0,u)\leq 2(\sharp S)e^{-\epsilon^2 C N}$$ }
for all $m_0\in \Delta_d,$ all every stationary strategy $u.$
\end{lem}
Since $C$ is independent of $N,$ and $(e^{-\epsilon^2 C})^N$ is summable, we can use the dominated convergence theorem:
     for all $\epsilon > 0$,

        { $$\sum_N\mathbb{P} \left(sup_{0\leq\tau\leq T} \parallel M^N(\tau)-m(\tau)\parallel_{\infty}>\epsilon |  \ {M}^{N}(0)=m_0,u \right) < \infty,$$ }
        By Borel-Cantelli's lemma,
for every fixed $t<\infty,$ the random variable $\nu^{N,t}:=\sup_{0\leq\tau\leq t} \parallel M^N(\tau)-m(\tau)\parallel_{\infty}$ converges almost completely towards $0.$ This $\nu^{N,t}$ implies that  converges almost surely to $0.$

We introduce of $\tilde X^N_j$ by scaling with step size $\frac{1}{N}.$ Then, $Z^N=(X^N,M^N)$ is approximate in some sense by a discrete time process $\tilde Z^N=(\tilde X^N,\tilde m^N)$ where $\tilde m^N(k)=m(\lfloor Nt\rfloor)$ $m$ solution of the ODE
        where $\tilde X_j^N$ is the discrete time jump process with transition matrix  $$q^N_{i,i+r}(\tilde m^N(k))=\frac{1}{N}L_{i,i+r}(m(\frac{k}{N}),u)).$$
Using the lemma~\ref{lemm} and
uniform Lipschitz continuity of of $L^N$, we obtain that
{
$$
\sup_{i,j} \sup_{0\leq\tau\leq t} \parallel q^N_{i,j}( M^N(\tau))-q_{i,j}(m(\tau))\parallel
$$
$$
\leq K(\epsilon_N+ \sup_{0\leq\tau\leq t}\|M^N(\tau)-m(\tau)\|).
$$ }
 Hence, we can write $\|M^N(\tau)-m(\tau)\| \leq K(\epsilon_N+\frac{1}{N^2})$ over set of event $\Omega_{\epsilon}=\{\|M^N(\tau)-m(\tau)\|\leq \epsilon \}$ and $ P(\Omega_{\epsilon})\geq 1-2(\sharp S)e^{-\epsilon^2 C N}\rightarrow 1.$
Thus, \bearn P({X^N}_{j,|[0,t]}={\tilde X^N}_{j,|[0,t]}| k \ \mbox{transitions})\geq \mathbb{E}(\epsilon^{Bin(\frac{1}{N},Nt)})
\\ \mathbb{E}(\epsilon^{Bin(\frac{1}{N},Nt)})=(1-\frac{1}{N}+\frac{1}{N}\epsilon)^{Nt}
\\ P({X^N}_{j,|[0,t]}={\tilde X^N}_{j,|[0,t]}| k \ \mbox{transitions})\geq e^{\epsilon} \eearn and this holds for any  $\epsilon$ arbitrary small.
We define $d(X,Y)=\sum_{k=0}\frac{1}{2^k}d(X_k,Y_k)$ where $d(X_k,Y_k)=\ind{X_k\neq Y_k}.$
Then, $d({X^N}_{j,|[0,t]},{\tilde X^N}_{j,|[0,t]})\longrightarrow 0$ when $N$ goes to infinity.

\noindent{\bf Convergence of the discrete time process}
To prove the weak convergence of $(\tilde X^N_j,\tilde M^N),$ we check the following steps:
\begin{itemize}
\item the discrete time empirical measures $\tilde M^N$ are tight (follows from Sznitman for finite states) and converges to a martingale problem. The limit $\tilde m$ is deterministic measure and is solution of ODE which has the unique solution $m$ (given $m_0,u$). Thus, $\tilde m=m.$
\item Conditionally to $\tilde M^N,$ $\tilde X_j^N$ converges to a martingale problem. The jump and drift process  $\tilde X$ with time dependent transition  is given by the limit of the marginal of $A^N(.|\tilde M^N,m_0,x_0,u).$ We derive the weak convergence of $(\tilde X^N_j,\tilde M^N)$ to $(\tilde X,\tilde m)$ where $\tilde m$ is deterministic and $\tilde X$ is random. For this we use the Theorem 17.25 and its discrete time approximation in Theorem 17.28 pages 344-347 in Kallenberg.
\end{itemize}

\subsection*{Sketch of Proof of Theorem \ref{theo-rew-cont}}
Since Skorohod's topology is induced by a metric,
  it is sufficient to show that whenever
$(X^N_j,m^N) \to (x,m)$ in
Skorohod's topology, we have:
 \begin{eqnarray}
  \lim_{N\longrightarrow \infty}\int_0^{\infty}e^{-\beta t}r^N(v,X^N_j(t),m^N(t)) dt \nonumber\\ \nonumber
   =
   \int_0^{\infty}e^{-\beta t}r(v,x(t),m(t)) dt
 \end{eqnarray}

 By
\cite{ethier-kurtz-05},
 page $1 1 7$, there is some
 sequence of increasing bijections $\lambda_n$: $[0,\infty) \to [0,\infty)$
  s.t. $$
 \frac{\lambda_n(t)-\lambda_n(s) }{t-s} \to 1
\mbox{ uniformly in } t \mbox{ and } s
$$
$$
\mbox{and }\parallel y_n(t)-y(\lambda_n(t))\parallel \to 0 \mbox{
uniformly in } t $$ over compact subsets of $[0,\infty).$
Fix $\eps > 0$,  arbitrary and consider

 { \begin{eqnarray}
 h^N &:=& |
 \int_0^{\infty}e^{-\beta t}r^N(X^N(t),v,m^N(t)) dt \nonumber\\ & -&
   \int_0^{\infty}e^{-\beta t}r(x(t),v,m(t)) dt | \nonumber \\
    &\leq &  \int_0^{\infty}e^{-\beta t} |r^N(x^N(t),v,m^N(t)) - r(x(t),v,m(t)) | dt \nonumber
 \end{eqnarray}
 }

 First let $K=
\sup_{x \in \calS,v, m \in \Delta} |r(x,v,m)| < {\infty}$ by
hypothesis, and pick some time $T$ large enough such
that $e^{-\beta T}K/\beta \leq \epsilon/3.$ Thus
 \be
 h^N \leq \eps/3 +
   \int_0^{T}e^{-\beta t}|
r(x^N(t),v,m^N(t)) - r(x(t),v,m(t))|
   dt \label{eq-ldskk}
 \ee

Second, we use the distance on $E$ defined by
 \be d((x,m),(x',m')) = \parallel m-m'\parallel+ \ind{x \neq x'}
 \ee
$$
\mbox{Let }
K' = \sup_{x \in \calS, v, m\in
\Delta_d}\frac{|r(x,v,m)-r(x',v,m')|}{\parallel m-m'\parallel}<
{\infty}
$$
by hypothesis. It is easy to see that for
all $x,x' \in \calS$ and $m, m' \in \Delta_d$:
 \be \parallel r(x,v,,m)-r(x',v,m')\parallel \leq K' d((x,m),(x',m'))
 \ee
  Thus, by \eref{eq-ldskk}:
  \be
  h^N \leq \eps/3 + K' \int_0^T e^{-\beta
  t}d\left((x^N(t),m^N(t)),(x(t),m(t))
  \right)dt
  \ee

  By \cite{ethier-kurtz-05},
 page $1 1 7$, there is some
 sequence of increasing bijections $\lambda^N$: $[0,\infty) \to [0,\infty)$
  s.t.
\[
 \frac{\lambda^N(t)-\lambda^N(s) }{t-s} \to 1
\mbox{ uniformly in } t \mbox{ and } s
\]
\[
\mbox{and }
  d\left((x^N(t), m^N(t),
 \; \;(x^N(\lambda^N(t)),m^N(\lambda^N(t)))\right)\to 0
\]
$ \mbox{
uniformly in } t \mbox{ over compact subsets of
}[0,\infty). $
 Thus there is some $N_0 \in \mathbb{N}$ such
that for $N \geq N_0$ and $t \in [0,T]$:
 \be
d\left((x^N(t), m^N(t),
 \; \;(x^N(\lambda^N(t)),m^N(\lambda^N(t)))\right)
 \leq \frac{\epsilon  \beta e^{\beta T}}{3 K'}
 \ee Thus,
by the triangular inequality for $d$: $h^N\leq $

 {
\[
 \leq \frac{ \epsilon}{3} + K'\int_0^T e^{-\beta
  t}d\left((x^N(t)m^N(t)),(x(\lambda^N(t)),m(\lambda^N(t))
  \right)dt
\]
\[
  + K' \int_0^T e^{-\beta
  t}d\left((x(\lambda^N(t)),m(\lambda^N(t))),(x(t),m(t)
  \right)dt
\]
\begin{equation}
  \leq
	\frac{ 2 \epsilon}{3} + K' \int_0^T e^{-\beta t} d
  \left(
        (x(\lambda^N(t)),m(\lambda^N(t))),(x(t),m(t))
  \right) \ dt
  \label{eq-sajjhfd}
\end{equation}
}

Third, let $D$ be the set of discontinuity points of $(x,m)$.
Since $(x,m)$ is cadlag, $D$ is enumerable, thus it is
negligible for the Lebesgue measure and
 \bearn
\lefteqn{\int_0^T e^{-\beta
  t}d\left((x(\lambda^N(t)),a,m(\lambda^N(t))),(x(t),a,m(t)
  \right)dt} \\
  &=& \int_0^T e^{-\beta
  t}d\left((x(\lambda^N(t)),m(\lambda^N(t))),(x(t),m(t)
  \right)  \ind{t \notin D} dt
\eearn

 Now $\lim_{N\longrightarrow\infty}\ \lambda^N(t)=t$
 and thus for $t \notin D$
 \ben \lim_{N\longrightarrow\infty}d\left((x(\lambda^N(t)),m(\lambda^N(t))),
 (x(t),m(t)
  \right) =0\een
  and thus by dominated convergence
 \be
 \lim_{N\longrightarrow\infty}\int_0^T e^{-\beta
  t}d\left((x(\lambda^N(t)),m(\lambda^N(t))),(x(t),m(t))
  \right)dt
  =
  0
 \ee
and for $N$ large enough the second term in the
right-hand side of \eref{eq-sajjhfd} can be made
smaller than $\eps/3$. Finally, for $N$ large enough,
$h^N \leq \eps$. This completes the proof.

\subsection*{Sketch of Proof of Theorem \ref{mainthm}}

Define the discounted stochastic evolutionary game with random number of interacting players in each local interaction in which each player in $x$ with the mixed action $u(.|x)$ receives $r(u,x, m(t))$ where $m(t)$ is the population profile at $t,$ which evolves under the dynamical system  (\ref{eq-ode}) and the between states follows the transition kernel $L.$ Then, a strategy of a player is the same as in the microscopic case and the discounted payoffs $$R(u_1,u_2,s_0,m_0)=\int_0^{\infty}e^{-\beta t}r(s(t),u_1,m[u_2](t)) dt$$ is the limit of $R^N(u_1,u_2,s_0,m_0)$ when $N$ goes to infinity, where $m[u_2]$ is the solution of the ODE $\dot{m}=f(u_2,m), m(0)=m_0$ . It follows that the asymptotic regime of the microscopic  game and the Markov decision evolutionary game (macroscopic game) are equivalent.

\subsection*{Sketch of Proof of Theorem \ref{prop3}}
We show that for every discount factor $\beta>0$ the optimal control problem $(OPT_N)$ (resp. the fixed-point problem $(FIX_N)$) has at least one $0-$optimal strategy.  It follows from the existence of equilibria in stationary strategies for finite stochastic games with discounted payoff:
The set of pure strategies is a compact space in the product topology (Tykhonov theorem). Thus, the set of behavioral strategies $\Sigma_j$ is a compact space and also convex as the set of probabilities on the pure strategies. For every player $j$ and every strategy profile $\sigma$ the marginal of the payoffs and constraints functions are continuous for any $\beta>0:$
$ \alpha_j\longmapsto R^N_j(\alpha_j,\sigma_{-j},s,m_0).$
Moreover, the stationary strategies is convex, compact and upper and lower hemi-continuous (as a correspondence).
Define $$\gamma_j(s,m_0,\sigma)=\arg\max_{\alpha_j\in \calU_s}\ R^N_{j}(\alpha_j,\sigma_{-j},s,m_0).$$ Then, $\gamma_j(m_0,\sigma) \subseteq\Sigma_j$ is a non-empty, convex and compact set and the product correspondence $$\gamma: \sigma\longmapsto (\gamma_1(s,m_{0},\sigma),\ldots,\gamma_N(s,m_{0},\sigma))$$ is upper hemi-continuous (its graph is closed). We now use the  Glicksberg generalization of Kakutani fixed point theorem, and there is a stationary strategy profile $\sigma^*$ such that $$\sigma^*\in \gamma(s,m_0,\sigma^*).$$
Moreover, if the game has  symmetric payoffs  and strategies for each type, there is a symmetric per type stationary equilibrium. This completes the proof.

\subsection*{Sketch of Proof of Theorem \ref{propt4}}
Let  $(U^N)_{N}$ be a sequence of solution of $(FIX_N)$ i.e equilibrium in the system with $N$ players.
Choose a subsequence $ N_k$ such that $U^{N_k}$ converges to some point $u$ when $k$ goes to infinity.
We can write $$ R^{N_k}(U^{N_k},U^{N_k})-R(U,U)= R^{N_k}(U^{N_k},U^{N_k})-R^{N_k}(U,U)+R^{N_k}(U,U)-R(U,U)
.$$ Since  $R^N(.,.)$  is continuous and converges uniformly to $R(.,.)$,  $R^{N_k}$ converges uniformly to $R,$ the second term $R^{N_k}(U,U)-R(U,U)\longrightarrow 0$ when $N_k \longrightarrow\infty$ and the first term $R^{N_k}(U^{N_k},U^{N_k})-R^{N_k}(U,U)$ can be  rewritten as
$ R^{N_k}(U^{N_k},U^{N_k})-R^{N_k}(U,U)= R^{N_k}(U^{N_k},U^{N_k})-R(U^{N_k},U^{N_k})
 + R(U^{N_k},U^{N_k})-R(U,U)+ R(U,U)-R^{N_k}(U,U).$ Each term goes to zero by continuity of $R,$ convergence of $U^{N_k}$ to $U,$ and by uniform convergence of $R^{N}$ to $R.$ Let $U^N$ be a $\epsilon_N-$equilibrium. Then, $ R^{N}(U^N,U^N)\geq R^N(v,U^N)-\epsilon_N,\ \forall v.$
Then any limit $U$ of a subsequence of $U^N$ satisfies $R(U,U)\geq R(v,U)-\epsilon,\ \forall v.$ Similarly, if $$ R^{N}(U^N,U^N)\geq R^N(v,v)-\epsilon_N,\ \forall v$$ then  any omega-limit $U$ of the sequence of $U^N$ satisfies $R(U,U)\geq R(v,v)-\epsilon,\ \forall v$ i.e $U$ is an $\epsilon-$optimal strategy. In particular if $(U^N)_{N}$ is a sequence of $\epsilon_N-$equilibria (resp. optimal strategies) with $\epsilon_N\longrightarrow 0$ when $N$ goes to infinity then any accumulation point $U$  of $(U^N)_{N}$ is a $0-$equilibrium (resp. $0-$optimal strategy).

      \end{document}